\documentclass{article-hermes}

\usepackage{graphicx}

\usepackage{amsmath}
\makeatletter
  \def\tagform@#1{\maketag@@@{[#1]\@@italiccorr}}
\makeatother

\usepackage[bookmarks=true, bookmarksnumbered=true, bookmarksopen=true,
		unicode=true, colorlinks=true,
		pagebackref=true,
		linkcolor=blue,citecolor=blue,filecolor=blue,urlcolor=blue
		]{hyperref}

\begin{document}

%\maketitle
\title[Une nouvelle mesure d'évaluation]%
      {Une nouvelle mesure pour l'évaluation des méthodes de détection de communautés}

\author{Vincent Labatut}
%\author{ }

% caractères turcs: http://ergut.wordpress.com/2008/03/16/hello-world/
\address{Département d'informatique, Université Galatasaray\\
Galatasaray Üniversitesi, Ç{\i}ra\u{g}an cad. n°36\\
Ortaköy 34357, \.{I}stanbul, Turquie\\[3pt]
vlabatut@gsu.edu.tr}
%\address{ }

\resume{La détection de communautés dans un réseau complexe est une tâche que l'on peut rapprocher de la classification non-supervisée réalisée en fouille de données classique. Pour cette raison, l'évaluation des algorithmes accomplissant ce type de traitement s'est faite jusqu'ici exclusivement au moyen de mesures comparables à celles utilisées en fouille de données. Cependant, dans le cas de l'analyse de réseau, celles-ci n'exploitent pas toute l'information disponible et sont susceptibles de fournir des résultats biaisés. Dans cet article, nous illustrons cette limitation et proposons une solution en modifiant une mesure existante. Nous l'appliquons ensuite à des données réalistes afin d'en effectuer une première évaluation expérimentale.}

\abstract{Community detection can be considered as a variant of cluster analysis applied to complex networks. For this reason, all existing studies have been using tools derived from this field when evaluating community detection algorithms. However, those are not completely relevant in the context of network analysis, because they ignore a part of the available information, and can therefore lead to incorrect interpretations. In this article, we illustrate this limitation, and propose a solution by modifying an existing measure. We then apply it to realistic community-structured networks, in order to perform a first evaluation.}

\motscles{Réseaux complexes, détection de communautés, mesures d'évaluation.}

\keywords{Complex networks, community detection, evaluation measure.}

%\submitted[25/06/2012]{MARAMI 2012}{1}
\proceedings{MARAMI 2012}{-}

\maketitlepage

\section{Introduction}
%\vspace*{-1.3cm} % tentative pour éviter l'espacement abusif >> échec !
La détection de communautés est un champ de l'analyse de réseaux complexes, consistant à en caractériser la structure à un niveau mésoscopique, i.e. qui n'est ni celui du n\oe{}ud (microscopique) ni celui du réseau entier (macroscopique). Concrètement, on cherche à décomposer le réseau en un ensemble de sous-graphes interconnectés, chacun constituant ce que l'on appelle une communauté. Le problème n'a pour l'instant pas été posé formellement, mais les auteurs s'accordent sur une définition intuitive de l'objectif de cette tâche, qui est d'obtenir des communautés dont les n\oe{}uds sont densément interconnectés, par rapport au reste du réseau \cite{Fortunato2010}.
 
Une étude présentant une nouvelle méthode de détection de communautés s'articule généralement de la manière suivante : après avoir décrit leur algorithme, les auteurs l'appliquent à des réseaux réels et/ou artificiels, afin de comparer ses performances à celles de certains algorithmes existants. Cette procédure soulève plusieurs problèmes méthodologiques abordés ailleurs, portant notamment sur la nature des réseaux utilisés lors des tests, la définition de structures de communautés de référence, ou la généralisation des résultats obtenus à d'autres réseaux. Dans cet article, on s'intéresse plutôt à l'outil utilisé pour quantifier la performance de l'algorithme. Il s'agit toujours d'une mesure permettant d'associer un score à la structure de communautés estimée, en la comparant à la structure de communautés réelle, qui est connue (soit qu'elle ait été déterminée manuellement dans le cas d'un réseau réel, soit par construction pour un réseau artificiel). Pour ce faire, chacune des deux structures de communautés est considérée comme une simple partition de l'ensemble des n\oe{}uds, et le rôle de la mesure est de quantifier leur similarité.

Cependant, comme l'a relevé une étude récente \cite{Orman2012}, cette approche a ses limites. En effet, en ignorant totalement l'information topologique lors de l'évaluation des performances, elle court le risque de mettre sur le même plan des résultats de qualités très contrastées. Ainsi, deux structures dont les communautés ont des propriétés topologiques significativement différentes peuvent aboutir à des scores similaires. Dans cet article, nous proposons de modifier la pureté, qui est l'une de ces mesures traditionnelles, afin de lui faire tenir compte de la distribution des liens du réseau. Notre but est d'obtenir un outil permettant une discrimination plus pertinente des algorithmes. Dans la section suivante, nous passons rapidement en revue l'approche utilisée traditionnellement dans la littérature pour évaluer les performances de ce type d'algorithme, et décrivons plus particulièrement la mesure de pureté. Nous soulignons ensuite la limite de cette approche dans le contexte qui est le nôtre, avant d'élaborer notre nouvelle mesure. Nous la validons alors expérimentalement, et concluons par une discussion de notre travail et de ses développements possibles.

\section{Approche traditionnelle}
\label{sec:revue}
Il existe en fouille de données classique une tâche appelée \emph{classification non-supervisée} (ou cluster analysis), consistant à partitionner un ensemble d'objets afin d'identifier des groupes homogènes. Chaque objet est décrit individuellement au moyen d'un ensemble d'attributs, et la procédure est menée en comparant les objets entre eux grâce à ces attributs. Le parallèle avec la détection de communautés est assez flagrant : dans le cas des réseaux complexes, les objets sont des n\oe{}uds, et l'information dont on dispose prend la forme de liens. La différence entre les deux tâches est donc qu'au lieu de considérer une information individuelle (les attributs), on exploite une information relationnelle (les liens). Cependant, dans les deux cas le résultat est le même : une partition de l'ensemble des objets, qui est appelée \emph{structure de communautés} dans le cadre de l'analyse de réseaux complexes.

Il n'est donc pas surprenant que les auteurs ayant développé des algorithmes de détection de communautés se soient tournés vers les méthodes utilisées en fouille de données classique pour évaluer les outils de partitionnement \cite{Danon2005}, ou bien aient abouti à la définition de méthodes similaires \cite{Girvan2002}. En classification, cette opération est réalisée grâce à une mesure permettant d'obtenir un score correspondant à la performance de l'outil. Lorsqu'une partition de référence est disponible, ce score représente la similarité entre cette partition réelle et celle estimée par l'algorithme considéré ; on parle alors de \emph{critère d'évaluation externe} \cite{Manning2008}. Cependant, il existe une multitude de mesures de ce type, et dans le domaine de la classification, le débat concernant le choix de l'outil d'évaluation le plus approprié est encore d'actualité \cite{Liu2007} : ceci montre bien l'importance de ce point méthodologique. En effet, à quoi bon évaluer un outil si la méthode utilisée n'est pas pertinente ? 

La plupart des mesures utilisées en classification ont également été appliquées en détection de communautés. Dans cet article, nous nous concentrons sur la pureté \cite{Manning2008}, qui est historiquement la première à avoir été utilisée dans ce contexte \cite{Girvan2002}, mais surtout la plus simple à adapter selon la méthode que nous développons plus loin. Considérons deux partitions $X=\{x_1,...,x_I \}$ et $Y=\{y_1,...,y_J \}$ d'un même ensemble, où les $x_i$ et $y_j$ sont les parties et $n$ est le nombre total d'éléments dans l'ensemble partitionné. La pureté d'une partie $x_i$ s'exprime par rapport à l'autre partition $Y$ de la manière suivante :

\begin{equation}
\label{f:PuretePartie}
Pur(x_i,Y)=\max_j \frac{| x_i \cap{} y_j |}{| x_i |}
\end{equation}

En d'autres termes, on cherche parmi toutes les parties de $Y$ laquelle est la plus présente dans $x_i$, et on calcule quelle proportion des éléments de $x_i$ elle représente. Plus l'intersection entre les deux parties est grande et plus la pureté est élevée, et leur correspondance forte. La pureté totale de la partition $X$ par rapport à la partition $Y$ est obtenue en effectuant la somme des puretés de chaque $x_i$, pondérée par leur prévalence dans l'ensemble traité :

\begin{equation}
\label{f:PureteTotale}
Pur(X,Y) = \sum_i \frac{| x_i |}{n}Pur(x_i,Y)
\end{equation}

Il est important de remarquer que la pureté n'est pas une mesure symétrique : calculer la pureté de $Y$ par rapport à $X$ revient à s'intéresser aux parties de $X$ majoritaires dans chaque partie de $Y$, et on n'a donc pas, sauf cas particulier, d'égalité entre $Pur(X,Y)$ et $Pur(Y,X)$.

Du point de vue de la détection de communautés, on peut donc utiliser deux mesures distinctes, suivant que l'on calcule la pureté des communautés estimées par rapport aux communautés réelles, ou le contraire. En classification, c'est généralement la première qui est appliquée, et appelée simplement \emph{pureté}, alors que la seconde est la \emph{pureté inverse} \cite{Manning2008}. Il est difficile de déterminer laquelle de ces deux versions a été effectivement utilisée en détection de communautés. En effet, dans leur article fondateur, Girvan et Newman donnent une définition très succincte de la mesure qu'ils utilisent, appelée \emph{fraction de n\oe{}uds correctement classifiés} \cite{Girvan2002}. Un article ultérieur semble indiquer qu'il s'agit de la pureté inverse \cite{Newman2004} (note 19), ce qui nous a été confirmé directement par Newman. Un grand nombre de travaux menés par la suite par d'autres auteurs ont utilisé des mesures portant ce même nom. Cependant, en raison de l'imprécision initiale, il est tout à fait possible qu'ils aient en réalité utilisé la pureté non-inverse. C'est ce qui ressort, par exemple, d'un commentaire que Danon \textit{et al}. \cite{Danon2005} (p.4) adressent à la mesure de Newman (donc à la pureté inverse), mais qui n'est en réalité valide que pour la pureté non-inverse.

Un autre point important est que chacune de ces deux versions de la mesure souffre d'une limitation. La pureté a tendance à favoriser les algorithmes identifiant de nombreuses communautés de petite taille. Dans le cas extrême, si l'algorithme identifie $n$ communautés contenant chacune un seul n\oe{}ud, on obtient une pureté maximale, car chaque communauté estimée est complètement pure. Au contraire, la pureté inverse favorise les algorithmes identifiant peu de communautés de grande taille. Cette fois, le cas extrême est celui où l'algorithme identifie une seule communauté contenant tous les n\oe{}uds. On obtient là encore une pureté maximale, car chaque communauté réelle est en effet complètement pure, puisque tous les n\oe{}uds qu'elle contient appartiennent à la même (unique) communauté estimée. Pour pallier ce défaut, Newman a introduit une contrainte supplémentaire : quand une communauté estimée est majoritaire dans plusieurs communautés réelles, tous les n\oe{}uds concernés sont considérés comme mal classifiés.

La méthode utilisée en classification consiste plutôt à calculer la F-mesure, qui est la moyenne harmonique des deux versions de la pureté \cite{Artiles2007} :

\begin{equation}
\label{f:FMesure}
F(X,Y) = 2 \frac{Pur(X,Y) \cdot Pur(Y,X)}{Pur(X,Y)+Pur(Y,X)}
\end{equation}

La mesure obtenue est symétrique, et cette combinaison est supposée compenser les biais précédemment décrits. Cette approche pénalise de façon similaire la sous-estimation et la surestimation du nombre de communautés, c'est pourquoi nous nous baserons sur elle par la suite, et non pas sur l'ajustement de Newman. 

\section{Limitation des mesures existantes}
Les mesures issues de la fouille de données classique, pureté comprise, ne considèrent par définition une structure de communauté que comme une partition de l'ensemble des n\oe{}uds. Ceci peut poser problème, car toutes les erreurs de classification n'ont pas forcément la même importance. Considérons par exemple la Figure 1, qui représente un réseau simple constitué de deux communautés, chacune représentée par une couleur différente. On propose deux estimations $A$ et $B$ de cette structure de communautés de référence, elle-même notée $R$. Chacune comporte une erreur de classification : un n\oe{}ud appartenant à la communauté de gauche est placé dans celle de droite. Pour $A$, il s'agit du n\oe{}ud 2 alors que pour $B$ c'est le n\oe{}ud 6. Si on applique les mesures présentées dans la section précédente, de manière à comparer $A$ et $R$, on obtient un score de $0,9$ à la fois pour la pureté et la pureté inverse. Par conséquent, la F-mesure, qui est leur moyenne, atteint la même valeur. Rien ne permettant de distinguer $A$ et $B$ en termes de partitions, on obtient bien entendu exactement les mêmes scores quand on compare $B$ et $R$ en utilisant les mêmes mesures. En d'autres termes, aucune de ces mesures ne permet de distinguer les deux erreurs représentées respectivement par $A$ et $B$.

Pourtant, intuitivement ces erreurs ne semblent pas du tout équivalentes. En effet, le n\oe{}ud 2 est bien plus ancré dans sa communauté réelle, et parait donc beaucoup plus important que le n\oe{}ud 6. L'erreur le concernant est donc plus grave, et par conséquent l'estimation $B$ est de meilleure qualité, et son score devrait être meilleur. Pour éprouver objectivement cette intuition, on peut s'intéresser à la \emph{modularité} \cite{Newman2004a} des trois partitions. Cette mesure quantifie la qualité d'une partition de façon aveugle, i.e. sans tenir compte d'une référence. Pour ce faire, elle compare la proportion de liens situés à l'intérieur des communautés, à l'espérance de cette même proportion calculée pour un modèle aléatoire générant des réseaux similaires (même taille et distribution de degré). La modularité est utilisée comme fonction objectif par de nombreux algorithmes de détection de communautés \cite{Fortunato2010}. Dans notre cas, la référence $R$ atteint une modularité de $0,36$, $A$ de $0,18$ et $B$ de $0,34$.  Plus que leur ordre de grandeur, ce qui nous intéresse ici sont les différences relatives entre ces valeurs : $A$ obtient une valeur bien inférieure à $B$, ce qui confirme notre intuition.
 
\begin{figure}[h]
{\setlength{\fboxsep}{2mm}
\fbox{\parbox{116mm}{%
\begin{center} \small 
\includegraphics[width=0.90\textwidth]{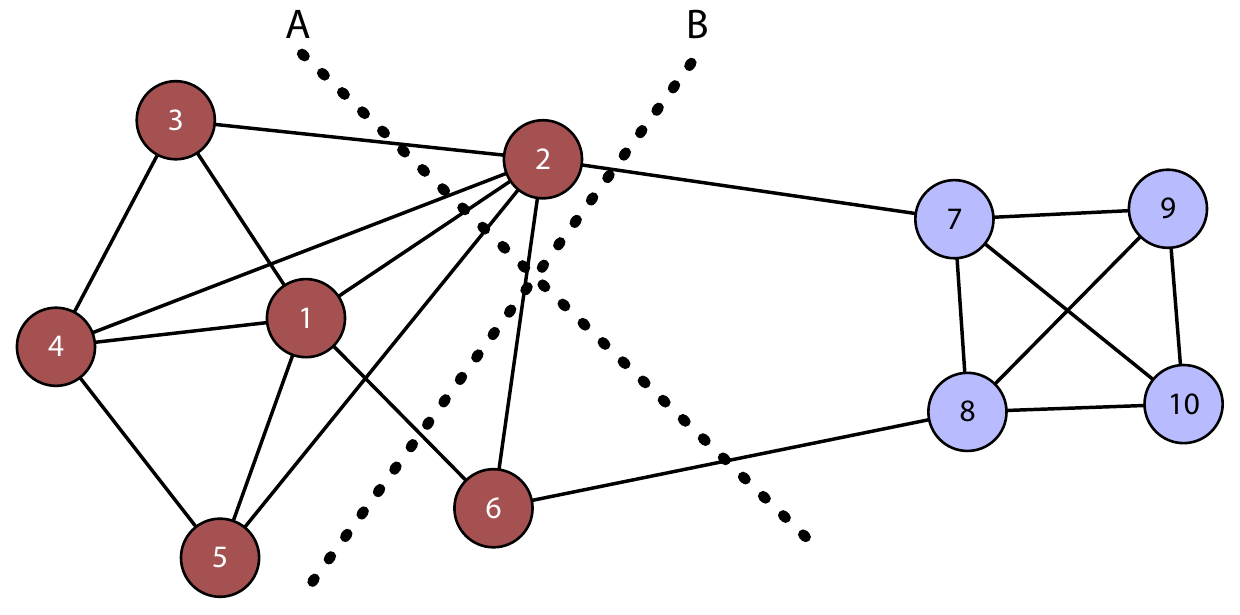}
\end{center}\vspace{-6pt}
}}
\caption{Cas illustrant la limitation de l'évaluation exclusivement basée sur les partitions. Les couleurs symbolisent les communautés réelles ; les droites $A$ et $B$ représentent les divisions correspondant à deux structures de communautés estimées différentes.}
\label{fig:Limitation}
}
\end{figure}

Cette observation effectuée sur un exemple simple est corroborée sur des réseaux plus réalistes, par une étude récente \cite{Orman2012}. Celle-ci compare des structures de communautés en considérant plusieurs propriétés topologiques permettant de les caractériser. Elle conclut qu'un algorithme de détection de communautés peut atteindre pour un réseau donné un score très élevé, sans que cela reflète forcément une similarité forte, d'un point de vue topologique, entre la structure de communautés estimée et la référence. Le corollaire est que deux algorithmes peuvent atteindre des scores très proches, mais correspondre à des structures de communautés topologiquement bien différentes. On peut en déduire que les mesures existantes ne sont pleinement adaptées ni à l'évaluation des performances d'un algorithme dans l'absolu, ni à la comparaison de plusieurs algorithmes.

\section{Proposition d'une nouvelle mesure}
Le problème vient évidement du fait que considérer qu'une structure de communautés est une simple partition revient à ignorer une partie de l'information disponible : la topologie du réseau. Pour effectuer une évaluation plus pertinente, les auteurs de l'étude mentionnée précédemment proposent d'utiliser de façon conjointe les mesures traditionnelles et différentes propriétés topologiques \cite{Orman2012}. Cependant, ils reconnaissent eux-mêmes que l'utilisation de ces dernières n'est pas aisée, car leur quantification prend la forme de plusieurs séries numériques, difficiles à comparer.

La solution que nous proposons ici est au contraire de modifier la pureté, de manière à ce qu'elle tienne compte de la topologie du réseau. On continue à bénéficier ainsi de la concision d'un score unique pour mesurer la performance de l'algorithme. Par rapport aux autres mesures existantes, la pureté présente l'avantage de pouvoir être reformulée de manière à faire apparaître individuellement la contribution de chaque n\oe{}ud à l'erreur totale. On définit d'abord la notion de pureté d'un n\oe{}ud $u$ pour une partition $X$ par rapport à une autre partition $Y$ :

\begin{equation}
\label{f:TopoPureteElement}
Pur(u,X,Y) = \delta{} (\operatorname*{arg\,max}_j |x_\alpha{} \cap{} y_j |, \beta{})
\end{equation}

Avec $u\in{}x_\alpha{}$ et $u\in{}y_\beta{}$  ; et où $\delta{}$ est le symbole de Kronecker (i.e. $\delta{}(a,b)=1$ si $a=b$, et $0$ sinon). La fonction est donc binaire : $1$ si $y_\beta{}$ est bien la partie de $Y$ majoritaire dans $x_\alpha{}$, et $0$ sinon. Considérons par exemple $Pur(2,A,R)$ dans le cas de la Figure \ref{fig:Limitation}, et notons $r_1$ la partie rouge de $R$ et $r_2$ sa partie bleue. Dans $R$, le n\oe{}ud 2 appartient à la partie rouge, on a donc $\beta=1$. Dans $A$, il appartient à la partie de droite, qui a une intersection plus grande avec $r_2$ qu'avec $r_1$, donc $Pur(2,A,R) = \delta{} (2,1) = 0$. Pour le n\oe{}ud 6, on aurait au contraire obtenu $Pur(6,A,R) = \delta{} (1,1) = 1$.

On peut ensuite reformuler la pureté d'une partie $x_i$ en moyennant la pureté de ses n\oe{}uds :

\begin{equation}
\label{f:TopoPuretePartie}
Pur(x_i,Y) = \frac{1}{|x_i|} \sum_{u\in{}x_i} Pur(u,X,Y) 
\end{equation}

Cette expression est équivalente à celle de l'équation [\ref{f:PuretePartie}], et permet donc de dériver la pureté totale de la partition comme en [\ref{f:PureteTotale}]. En développant, on obtient :

\begin{equation}
\label{f:TopoPureteTotale}
Pur(X,Y) = \sum_i \sum_{u\in{}x_i} \frac{1}{n} Pur(u,X,Y)
\end{equation}

On peut remarquer que la pureté de chaque n\oe{}ud est pondérée par la valeur $1/n$. On se propose d'introduire la prise en compte de l'aspect topologique dans cette mesure, en remplaçant ce poids uniforme par une valeur $w_u$ propre à chaque n\oe{}ud $u$, de manière à pénaliser plus fortement les erreurs concernant les n\oe{}uds les plus importants topologiquement :

\begin{equation}
\label{f:TopoPureteTotalePonderee}
Pur'(X,Y) = \sum_i \sum_{u\in{}x_i} \frac{w_u}{\sum\limits_v w_v} Pur(u,X,Y)
\end{equation}

La normalisation des poids par leur somme totale permet de garder une valeur finale comprise entre $0$ et $1$. Enfin, en appliquant à [\ref{f:TopoPureteTotalePonderee}] le même principe que dans [\ref{f:FMesure}], on peut définir une version modifiée de la F-mesure, prenant en compte la topologie du réseau, et que l'on note $F'$.

Toute la question est alors de savoir comment caractériser l'importance des n\oe{}uds. L'idée est qu'une erreur de classification concernant un n\oe{}ud fortement intégré à sa communauté doit compter plus que s'il s'agit d'un n\oe{}ud situé à la lisière. On peut pour cela considérer le degré $d$ du n\oe{}ud, en le divisant par le degré maximal du réseau pour le normaliser : $d(u) / \max_v d(v)$ . On donne ainsi plus de poids aux hubs communautaires tels que le n\oe{}ud 2 de la Figure \ref{fig:Limitation}, et moins aux n\oe{}uds plus périphériques comme le n\oe{}ud 6.

Cependant, on peut faire deux critiques à cette approche. Premièrement, il est possible qu'un n\oe{}ud de degré élevé ait des connexions distribuées sur de nombreuses communautés, ce qui empêche toute intégration forte avec l'une d'elle en particulier. Puisque l'appartenance de ce n\oe{}ud semble assez incertaine, lui donner un grand poids semble peu approprié. De plus, utiliser seulement le degré amène à minimiser l'importance des n\oe{}uds de degré très faible, mais dont l'entièreté des connexions est réalisée dans leur communauté. La mesure d'\emph{enchâssement} (ou embeddedness) \cite{Lancichinetti2010} permet de prendre ces deux traits en compte. Il s'agit du rapport entre le degré interne d'un n\oe{}ud $d_{int} (u)$ et son degré total $d(u)$. Le degré interne d'un n\oe{}ud correspond au nombre de ses voisins directs situés dans la même communauté que lui. La mesure est comprise entre $1$ (tous les liens du n\oe{}ud sont dans sa communauté) et $0$ (aucun lien dans sa communauté).

Afin de combiner les deux aspects, on se propose de multiplier ces deux valeurs (degré normalisé et enchâssement). Ainsi, plus un n\oe{}ud cumule ces deux propriétés et plus il est important pour nous. Après simplification, le poids d'un n\oe{}ud $u$ est donc défini de la façon suivante :

\begin{equation}
\label{f:Poids}
w_u = \frac{d_{int} (u)}{\max\limits_v d(v)} 
\end{equation}

Si on calcule $F'$ pour l'exemple de la Figure \ref{fig:Limitation}, on obtient les valeurs $0,85$ pour $A$ et $0,94$ pour $B$, ce qui correspond au comportement attendu. Bien sûr, la formule [\ref{f:TopoPureteTotalePonderee}] est générale et peut être adaptée en fonction des besoins, en modifiant l'expression du poids. Par exemple, si les liens sont pondérés, on peut considérer la force des n\oe{}uds (somme des poids des liens attachés au n\oe{}ud) plutôt que leur degré.

\section{Évaluation expérimentale}
Afin d'éprouver notre mesure sur des données plus réalistes, nous l'avons appliquée à des réseaux artificiels possédant une structure de communautés. La méthode et les paramètres de génération retenus sont exactement les mêmes que dans l'étude déjà citée \cite{Orman2012}, afin de bénéficier des observations qu'elle contient quant aux propriétés topologiques des structures de communautés réelles et estimées. Pour la même raison, nous avons appliqué les mêmes algorithmes de détection de communautés : Copra \cite{Gregory2010}, FastGreedy \cite{Newman2004}, InfoMap \cite{Rosvall2008}, InfoMod \cite{Rosvall2007}, Louvain \cite{Blondel2008}, MarkovCluster \cite{Dongen2008}, Oslom \cite{Lancichinetti2011} et WalkTrap \cite{Pons2005}. Notre but n'était pas d'identifier les meilleurs algorithmes, mais bien d'éprouver le pouvoir discriminant de notre mesure. Nos données sont constituées de $5$ réseaux de $n=25000$ n\oe{}uds, dont les principales propriétés topologiques sont compatibles avec les descriptions de réseaux réels de la littérature : distribution du degré, transitivité, tailles des communautés, enchâssement des n\oe{}uds, etc. (cf. \cite{Orman2012} pour plus de détails sur les algorithmes ou la méthode générative).

Les auteurs de l'étude précédemment citée utilisent un panel représentatif de mesures pour évaluer les performances des algorithmes. Cependant la fraction de n\oe{}uds correctement classifiés qu'ils ont retenue est la version de Newman, que nous avons décrite dans la section \ref{sec:revue}. Or, notre propre mesure est basée sur la F-mesure appliquée à la pureté classique. Afin que la comparaison entre notre mesure et la mesure non-topologique soit pertinente, nous avons donc recalculé les performances des algorithmes en utilisant la F-mesure non-topologique de l'équation [\ref{f:FMesure}], et bien sûr notre propre variante, la version topologique basée sur la formule [\ref{f:TopoPureteTotalePonderee}].

\begin{figure}[h]
{\setlength{\fboxsep}{2mm}
\fbox{\parbox{116mm}{%
\begin{center} \small 
\includegraphics[width=0.90\textwidth]{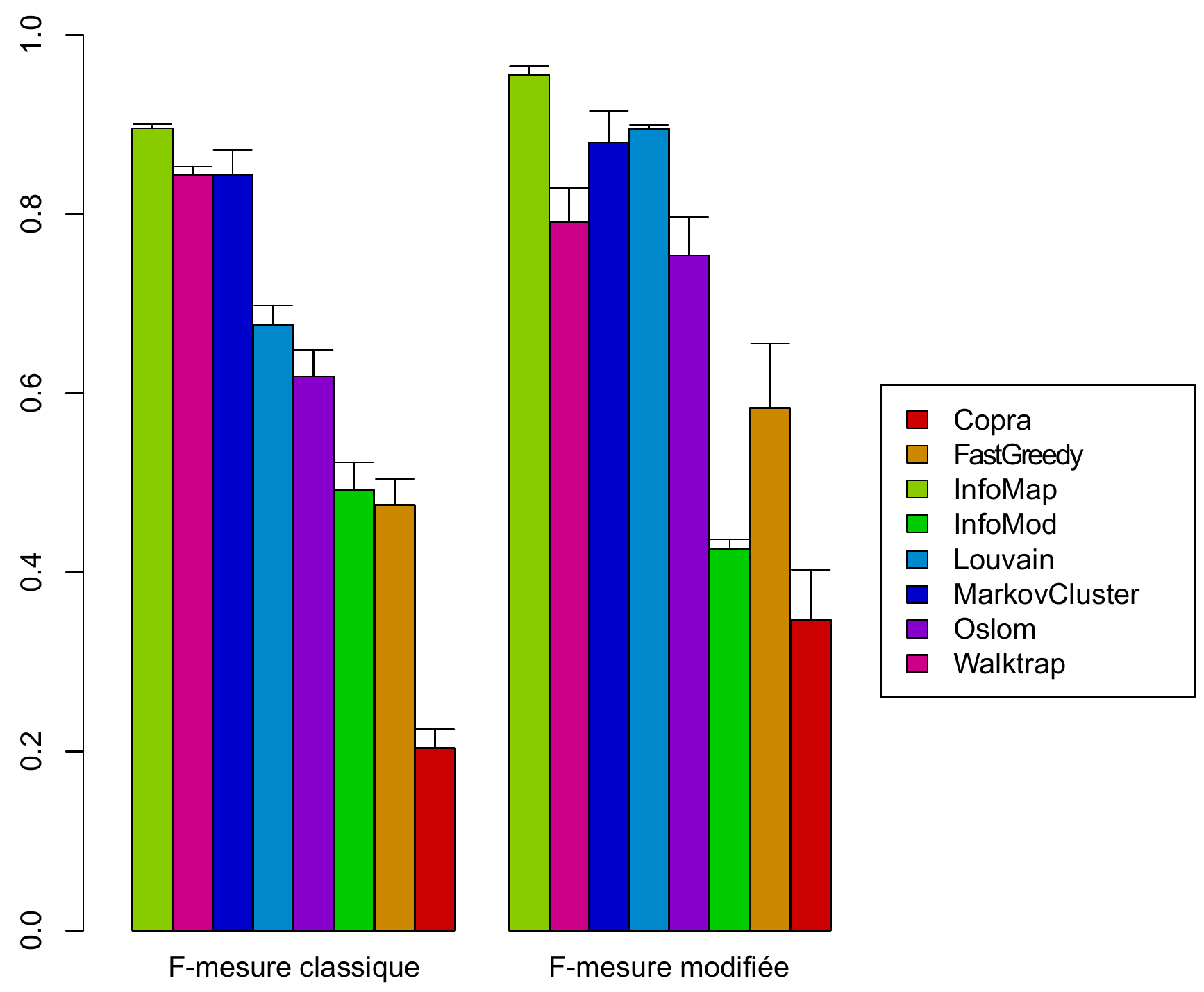}
\end{center}\vspace{-6pt}
}}
\caption{Comparaison des résultats obtenus avec la F-mesure classique et notre version modifiée pour tenir compte de la topologie.}
\label{fig:Resultats}
}
\end{figure}

La Figure \ref{fig:Resultats} présente les résultats obtenus. Seuls deux algorithmes, InfoMod et WalkTrap, voient leur score diminuer avec notre mesure, ce qui signifie que leurs performances sont majoritairement dues à des n\oe{}uds de poids faible, et donc considérés comme peu importants topologiquement parlant. On observe l'effet inverse pour tous les autres algorithmes.

Les classements des algorithmes selon les deux mesures, obtenus après Anova et test de Tukey avec un seuil de signification de $\alpha=0.05$, sont présentés dans le Tableau \ref{tab:Resultats}. Les algorithmes dont les résultats ne sont pas considérés comme significativement différents sont placés sur la même ligne. On peut voir que l'ordre ne change pas radicalement pour notre mesure. En particulier, les algorithmes ayant obtenu les scores les plus extrêmes restent les mêmes (InfoMap et Copra). Ceci est cohérent avec l'analyse effectuée dans \cite{Orman2012}, qui remarquait que les structures de communautés détectées par InfoMap sont extrêmement proches, pour l'ensemble des propriétés étudiées, de celles des structures réelles (et inversement pour Copra).

\begin{table}[h]
	\centering
	\begin{tabular}{|l|l|l|l|}
		\hline
 		\multicolumn{2}{|l|}{\textbf{F-mesure classique}} & 
 		\multicolumn{2}{|l|}{\textbf{F-mesure modifiée}} 	\\
		\hline
		Rang	&	Algorithmes		& Rang	&	Algorithmes					\\
		\hline
		1		&	InfoMap			& 1		&	InfoMap						\\
		\hline
		2		&	MarkovCluster, Walktrap	& 2	&	Louvain, MarkovCluster	\\
		\hline
		4		&	Louvain, Oslom	& 4		&	WalkTrap, Oslom				\\
		\hline
		6		&	InfoMod, FastGreedy	& 6		&	FastGreedy				\\
		\hline
		8		&	Copra			& 7		&	InfoMod, Copra				\\
		\hline

	\end{tabular}
	\caption{Classement des algorithmes en fonction de la F-mesure classique et de sa variante topologique.}
	\label{tab:Resultats}
\end{table}

InfoMod est rétrogradé au niveau de Copra, alors qu'il obtient des performances comparables à celles de FastGreedy avec la F-mesure classique. De la même façon, WalkTrap est rétrogradé au niveau d'Oslom alors qu'il était auparavant équivalent à MarkovCluster. Ce résultat est lui aussi cohérent avec les remarques faites dans \cite{Orman2012}, qui relèvent que WalkTrap est considéré comme très proche d'InfoMap au vu des mesures classiques, alors que les propriétés topologiques des communautés qu'il détecte sont plus éloignées de la référence, notamment en ce qui concerne la taille des communautés. La même remarque avait été faite au sujet de MarkovCluster, cependant notre mesure ne semble pas la confirmer. Tout au plus observe-t-on un accroissement de l'écart avec InfoMap, mais celui-ci était déjà significatif avec la mesure classique. Enfin, la performance mesurée pour Louvain augmente suffisamment pour le placer au niveau de MarkovCluster, alors qu'il était comparable à Oslom avec la mesure classique. Pour résumer, on peut conclure que, d'après ces premières évaluations, la mesure que nous proposons est en accord avec l'analyse précédente concernant la relation entre performance classique et propriétés topologiques des communautés estimées. Cependant, des manipulations supplémentaires sont nécessaires pour évaluer au mieux l'opportunité de la préférer aux mesures classiques.

\section{Conclusion}
Dans cet article, nous nous sommes intéressés aux mesures utilisées pour évaluer les algorithmes de détection de communautés. Toutes celles mentionnées dans la littérature sont du type de celles utilisées en fouille de données, et plus précisément en classification non-supervisée. Notre première contribution est d'avoir montré qu'aucune n'est pleinement appropriée à cette tâche, car elles ignorent complètement la topologie du réseau traité. Ceci diminue leur pertinence, et peut amener l'utilisateur à une interprétation erronée des scores obtenus. Notre deuxième contribution est d'avoir défini une variante de la pureté, la plus répandue de ces mesures, afin de résoudre ce problème. Nous avons pour cela introduit une modification générale consistant à y intégrer un poids pour chaque n\oe{}ud, ce qui permet de pénaliser individuellement chaque erreur de classification. L'adaptation à l'analyse des réseaux complexes est alors effectuée en utilisant une propriété topologique nodale comme poids. Dans un premier temps, nous avons proposé de tenir compte à la fois du degré du n\oe{}ud et de son enchâssement dans sa communauté réelle. Notre troisième contribution est d'avoir illustré l'apport de notre mesure, en l'utilisant sur les résultats fournis par une sélection d'algorithmes de détection de communautés, appliqués à des réseaux réalistes générés artificiellement. Le classement obtenu pour la mesure modifiée est cohérent avec une analyse approfondie de la topologie des structures de communautés estimées, qui avait été conduite précédemment \cite{Orman2012}. En d'autres termes, les performances relatives des algorithmes nous semblent, sur ces données, rendues de façon plus pertinente qu'avec la mesure classique.

L'une des limitations de ce travail concerne le choix de la définition des poids que nous avons introduits pour tenir compte de la topologie du réseau. En effet, chaque poids est supposé représenter l'importance du n\oe{}ud associé dans le réseau, et cette notion est difficile à formuler objectivement. Dans cet article, nous avons pénalisé les algorithmes incapable de traiter correctement les n\oe{}uds supposés faciles à classifier : ceux situés au c\oe{}ur des communautés. Mais il serait possible de prendre le contre-pied de cette approche, en faisant l'hypothèse que les algorithmes classifient tous correctement ce type de n\oe{}uds, et qu'il faudrait donner plus d'importance à ceux situés en lisière des communautés. On aboutirait alors probablement à des résultats assez différents. Une autre critique peut être formulée, mais à propos de la pureté cette fois. Dans le domaine de la fouille de données classique, certains auteurs préconisent d'appliquer une correction à ce type de mesures, afin de tenir compte de la part de chance dont le classificateur bénéficie \cite{Guggenmoos1993}. L'idée est de déterminer sa performance réelle, en déduisant du score obtenu les résultats corrects dus au seul hasard. Dans un premier temps, nous avons décidé de ne pas effectuer cette correction sur notre mesure, bien que cela soit tout à fait possible.

En plus de ces points, notre travail peut être étendu de multiples manières. On peut tout d'abord envisager deux autres types d'utilisation de notre mesure. Bien que celle-ci ait été définie pour comparer une partition estimée à une partition de référence, on pourrait l'utiliser pour comparer deux partitions estimées, comme il est déjà possible de le faire avec les mesures classiques. On effectuerait alors une comparaison qui tiendrait compte de l'aspect topologique. Elle pourrait également être appliquée dans le cadre de la classification classique, i.e. sur des données non-relationnelles, dans le cas où l'on jugerait nécessaire d'assigner différentes importances aux objets classifiés. La mesure pourrait aussi être modifiée afin de l'adapter aux particularités de réseaux spécifiques : liens orientés, liens pondérés, etc. Enfin, le principe même d'intégrer une prise en compte de l'information relationnelle dans une mesure classique pourrait être appliqué à d'autres mesures telles que l'\emph{information mutuelle normalisée} \cite{Strehl2002} ou l'\emph{index de Rand} \cite{Rand1971}, qui ont sur la pureté l'avantage d'être symétriques.

\nocite{*}
\bibliography{marami12}

\end{document}